\documentstyle[11pt]{article}
\def\eps{\varepsilon}
\begin{document}
\begin{center}
{\large\bf Cohomology of a Quaternionic Complex}\\[5mm]
Robin Horan\\
{\em School of Mathematics and Statistics}\\
{\em University of Plymouth}\\
{\em Plymouth, PL4 8AA, U.K,}\\
{\em e-mail: r.horan@plymouth.ac.uk}
\end{center}
\vspace{1cm}
\begin{abstract}
We investigate the cohomology of a certain elliptic complex defined on a
compact quaternionic-K\"{a}hler manifold with negative scalar curvature.
We show that this particular complex is exact, with the possible exception
of
one term.
\end{abstract}
Let $(M,g)$ be an oriented $4k$-dimensional compact
quaternionic-K\"{a}hler manifold
having negative scalar curvature. In \cite{hor} we proved a rigidity
theorem
{}for such manifolds which was, in essence, a consequence of the
vanishing of a
certain cohomolgy group on the twistor space $Z$, of $M$. The paper was
largely devoted to the proof of a vanishing theorem for
 the cohomolgy of the first term
of a particular complex on $M$, the required vanishing theorem on $Z$
being deduced
by the Penrose transform. This note is an extension of that work in
that we use
the same techniques to show that all the cohomology of the complex
vanishes, with the
possible exception of one term, so that the complex is exact except
possibly at the
right. We have not been able to establish whether the final term of
the complex
has non-trivial cohomology.

\noindent As in \cite{hor}, we shall make extensive use of the methods and
techniques developed by Bailey and Eastwoood in \cite{pcf}, these being
based on the abstract index notation of Roger Penrose \cite{sst}.
In this notation the indices are used
as `place-markers' and do not indicate a choice of basis. Thus the $E,H$
bundles of Salamon \cite{sal1}, are written as
${\cal O}^{^{A}}, {\cal O}^{^{A'}}$ respectively.
The (complexified) tangent bundle of $M$ is then ${\cal O}^{a} =
{\cal O} ^{^{A}}
\otimes {\cal O} ^{^{A'}} \stackrel{def}{=} {\cal O}^{^{AA'}}$.
The Levi-Civita
 connection on $M$ is written $\nabla _{a} =\nabla_{_{AA'}}$.
We shall take the Riemann curvature tensor to be given by
\begin{equation}
\label{riem}
2 \nabla_{[a} \nabla_{b]}w_{c}=R_{abc}^{\; \; \; \; \; d}w_{d}
\end{equation}
so that, in our convention, the standard metric on the sphere has
positive
curvature. (This differs from the convention adopted in \cite{pcf}
which was based on that of \cite{sst}. This latter is essentially a
relativity
book.)
The square brackets in the definition of the Riemann tensor, denotes
anti-symmetrisation. Round brackets will be used for symmetrisation.
 For more details on the notation and techniques used we refer the
reader to
\cite{sst}.

\noindent When written in abstract index notation, the complex we
consider
is defined on $M$ by the following

\begin{eqnarray}
&0\longrightarrow \Gamma(M,{\cal O}^{^{(A'_1 \ldots A'_n)}}) \stackrel
{\nabla_{B_{1}}^{A'_{n+1}}}{\longrightarrow}
\Gamma(M,{\cal O}_{_{B_1}}^{^{(A'_1 \ldots A'_{n+1})}})
\stackrel{\nabla_{B_{2}}^{A'_{n+2}}} \longrightarrow \Gamma(M,{\cal O}
_{_{[B_{1}B_{2}]}}^{^{(A'_1 \ldots A'_{n+2})}}) \longrightarrow \ldots&
\nonumber \\
&\ldots \stackrel{ \nabla_{B_{2k}}^{A'_{n+2k}}} \longrightarrow\Gamma(M,
{\cal O}_{_{[B_1 \ldots B_{2k}]}}^{^{(A'_1 \ldots A'_{n+2k})}})
\longrightarrow 0&
\label{cx1}
\end{eqnarray}
Here, for example, the bundle ${\cal O}_{_{[B_1B_2]}}^{^{(A'_1 \ldots
A'_{n+2})}}$ is
$\bigwedge ^{2} E^{*} \otimes S^{n+2}H$ (or its sheaf of smooth sections),
and $\Gamma (M, V)$ is the space of global smooth sections of the bundle
$V$ over $M$. The map $\nabla_{_{A}}^{^{A'_1}}$  is then $\varepsilon
^{^{A'_1A'}} \nabla _{a}$ with $\varepsilon ^{^{A'_1 A'}}$ being the
(covariantly constant) symplectic form on $H^{*}$.
This complex was studied by Salamon in \cite{sal2} and shown by him to
be an
elliptic complex. Since $M$ is compact we may use Hodge theory
\cite{wells}
to examine the cohomology of this complex.
It is well known that the cohomology of the first term vanishes for
negative
scalar curvature. The vanishing of the cohomology of the second term
was proved
in \cite{hor}, but can easily be obtained from the proof of
theorem~(\ref{thm1})
 with the appropriate changes.
For terms other than the last Hodge theory implies that
each cohomology class has a unique representative $f$ which satisfies
\begin{equation}
\begin{array}{ll}
\mbox{(a)} &  f \in
\Gamma(M,{\cal O}_{_{[B_1 \ldots B_m]}}^{^{(A'_1 \ldots A'_p)}})\\
\mbox{(b)} & \nabla_{_{[B_0}}^{^{(A'_0}}f_{_{B_1 \ldots B_m]}}
^{^{A'_1 \ldots A'_p)}}=0\\
\mbox{(c)} & \nabla_{_{AC'}}f^{^{AC'A'_2 \ldots A'_p}}
_{_{B_2\ldots B_m}}=0
\label{hodge1}
\end{array}
\end{equation}
where $1 < m < 2k$ and $p-m=n$.
If we let $D_0 = \nabla_{_{B_1}}^{^{A'_{n+1}}}$ and $D_1 =
\nabla_{_{B_2}}
^{^{A'_{n+2}}}$
 in (\ref{cx1}), then with $m=1$, (b) is the abstract index version of
$D_{1}f=0$ and (c) is that of $D_{0}^{\ast} f =0$,  $D_{0}^{\ast}$
being the
adjoint of the mapping $D_{0}$.

We shall take $u^{A}\overline{v}_{A}$ to be the positive definite inner-
product, where $\overline{v}$ is the quaternionic conjugate of $v$, on
the
bundle ${\cal O}^{A}$, so that $u^{A}\overline{u}_{A}$ is positive if
$u \neq 0$, and similarly for other tensors. The $L^{2}$-norm of a tensor
$v=v^{AB^{'}}_{\;\;\;\;CD^{'}}$ for example, is then
\begin{equation}
\parallel v \parallel^{2} = \int v^{AB^{'}CD^{'}}
\overline{v}_{AB^{'}CD^{'}}
\label{norm}
\end{equation}
with the integral taken over the manifold $M$. With this convention, one
can use integration by parts over $M$ and quickly obtain (c) above as the
abstract index version of $D_{0}^{\ast} f=0$, when $m=1$.
We shall show that these terms all vanish,
so that the above complex (\ref{cx1}) is exact, except possibly at the
right.
We note the following elementary facts.

\newtheorem{lemma1}{Lemma}[section]
\begin{lemma1}
\label{lemma1}
 If $S^{^{A'_0 \ldots A'_p}}$ is symmetric in its final $p$ indices
and $U_{_{B_0 \ldots B_m}}$ is antisymmetric in its final $m$ indices,
then
\begin{eqnarray}
 (p+1)S^{^{(A'_0 \ldots A'_p)}}&=&(p+1)S^{^{A'_0 \ldots
A'_p}} - \sum _{i=1}^{i=p}\eps^{^{A'_0A'_i}}S_{_{C'}}
^{^{\;\;\;C'A'_1 \ldots \hat{A'}_i \ldots A'_p}}
             \label{l1a}\\
(m+1)U_{_{[B_0 \ldots B_m]}}&=&U_{_{B_0 \ldots B_m}}+\sum_{i=1}^{i=m}
(-1)^iU_{_{B_iB_0 \ldots \hat{B}_i \ldots B_m}}
             \label{l1b}
\end{eqnarray}
\end{lemma1}

\noindent It is then quite easy to obtain the following.
\newtheorem{lemma2}[lemma1]{Lemma}
\begin{lemma2}
\label{lemma2}
If $f$ satisfies the conditions of (\ref{hodge1}) then
\begin{displaymath}
(p+1) \parallel \nabla _{_{[B_0}}^{^{A'_0}} f_{_{B_1 \ldots B_m]}}
^{^{A'_1 \ldots A'_p}} \parallel^2=p \parallel \nabla
_{_{ A'[B_0 }}
f_{_{B_1 \ldots B_m]}}^{^{A'A'_2 \ldots A'_p}} \parallel^2
\end{displaymath}

and $ \nabla ^{^{BA'_0}}f^{^{A'_1 \ldots A'_p}}_{_{BB_2 \ldots
B_m }}$ is symmetric in its primed indices.
\end{lemma2}
\noindent
{\em Proof} \quad For the first part we put $S^{^{A'_0 \ldots A'_p}}=
\nabla_{_{[B_0}}^{^{A'_0}}f_{_{B_1 \ldots B_m]}}
^{^{A'_1 \ldots A'_p}}$ in (\ref{l1a}) and using (b) of (\ref{hodge1})
 we get
\begin{displaymath}
(p+1)\nabla_{_{[B_0}}^{^{A'_0}}f_{_{B_1 \ldots B_m]}}
^{^{A'_1 \ldots A'_p}} = \sum_{i=1}^{i=p} \eps ^{^{A'_0A'_i}}
\nabla_{_{A'[B_0}}f_{_{B_1 \ldots B_m]}}^{^{A'A'_2 \ldots A'_p}}
\end{displaymath}
Contracting both sides of the above with
$\nabla^{^{[B_0}}_{_{A'_0}} \bar{f}^{^{\;B_1 \ldots B_m]}}
_{_{A'_1 \ldots A'_p}}$ and integrating the result over $M$,
one quickly obtains the first part. The second part can be
obtained by putting
$S^{^{A'_0 \ldots A'_p}}= \nabla ^{^{BA'_0}}f^{^{A'_1 \ldots
A'_p}}_{_{BB_2 \ldots B_m }}$ in (\ref{l1a}) and using (b) of
(\ref{hodge1}).

\noindent We can now state our main result.
\newtheorem{thm1}[lemma1]{Theorem}
\begin{thm1}
\label{thm1}
Let $M$ be a compact quaternionic-K\"{a}hler manifold of dimension $4k$,
for
$k>1$, and let $\Lambda =R/8k(k+2)$, where $R$ is its (non-zero) scalar
curvature. If $f$ is a harmonic element of
 $\Gamma(M,{\cal O}_{_{[B_1 \ldots B_m]}}^{^{(A'_1 \ldots A'_p)}})$,
 i.e. satisfies the conditions of(\ref{hodge1}),  and if $1<m<2k$, then
\begin{eqnarray}
\frac{(p+2)}{2(p+1)} \parallel \nabla_{_{A'}}^{^{[B_0}}f^{^{B_1
\ldots B_m]
A'A'_2 \ldots A'_p}}
\parallel^2& + & \frac{m}{2(m+1)} \parallel \nabla_{_{B}}^{^{C'}}
f^{^{BB_2 \ldots B_m A'_1 \ldots A'_p}} \parallel^2  \nonumber \\
& =& \Lambda \frac{(p+2)}{(m+1)}(2k-m)
\parallel f \parallel^2
\end{eqnarray}
Thus if $R<0$ then $f=0$.
\end{thm1}
\noindent {\em Proof}\quad
We have
\begin{eqnarray}
 \parallel \nabla_{_{A'}}^{^{[B_0}}f^{^{B_1 \ldots B_m]A'A'_2
\ldots A'_p}}
\parallel^2&=&\int ( \nabla_{_{K'}}^{^{[B_0}}f^{^{B_1 \ldots B_m]K'A'_2
\ldots
 A'_p}}) ( \nabla_{_{C'[B_0}} \bar{f}_{_{B_1 \ldots B_m]A'_2 \ldots A'_p}}
^{^{\;C'}} \nonumber \\
&=&-\int f^{^{B_1 \ldots B_m K'A'_2 \ldots A'_p}} \nabla^{^{B_0}}_{_{K'}}
 \nabla_{_{C'[B_0}} \bar{f}_{_{B_1 \ldots B_m]A'_2 \ldots A'_p}}^{^{\;C'}}
\nonumber \\
&=&-\int f^{^{B_1 \ldots B_m K'A'_2 \ldots A'_p}}( \nabla^{^{B_0}}_{_{
(K'}}
 \nabla_{_{C')[B_0}} \nonumber \\
&&+  \nabla^{^{B_0}}_{_{[K'}} \nabla_{_{C'][B_0}} ) \label{thm1a}
 \bar{f}_{_{B_1 \ldots B_m]A'_2 \ldots A'_p}}^{^{\;C'}}
\label{thm1b}
\end{eqnarray}
where the second line above is obtained from the first by using
integration
 by parts over the manifold. The latter integral on the right is
\begin{eqnarray*}
\frac{1}{2} \int f^{^{B_1 \ldots B_m K'A'_2 \ldots A'_p}} \eps_{_{K'C'}}
\nabla^{^{B_0}}_{_{D'}}
 \nabla_{_{[B_0}}^{^{D'}} \bar{f}_{_{B_1 \ldots B_m]A'_2 \ldots A'_p}}
^{^{\;C'}}&=&\frac{1}{2}\int f^{^{B_1 \ldots B_m A'_2 \ldots A'_p}}_{_{
C'}}
\nabla^{^{B_0}}_{_{D'}} \nabla_{_{[B_0}}^{^{D'}} \bar{f}_{_{B_1 \ldots
 B_m]A'
_2 \ldots A'_p}}^{^{\;C'}} \\
&=&-\frac{1}{2}  \parallel \nabla _{_{[B_0}}^{^{A'_0}} f_{_{B_1
\ldots B_m]}}
^{^{A'_1 \ldots A'_p}} \parallel^2
\end{eqnarray*}
the final expression being obtained from the previous one by using
integration
by parts again. Using the first part of lemma~(\ref{lemma2})
 and rearranging, (\ref{thm1b}) becomes
\begin{equation}
\frac{(p+2)}{2(p+1)} \parallel \nabla_{_{A'}}^{^{[B_0}}f^{^{B_1
\ldots B_m]A'
A'_2 \ldots A'_p}}\parallel^2=-\int f^{^{B_1 \ldots B_m K'A'_2
\ldots A'_p}}
 \nabla^{^{B_0}}_{_{(K'}}
 \nabla_{_{C')[B_0}} \bar{f}_{_{B_1 \ldots B_m]A'_2 \ldots A'_p}}^{^{
\;C'}}
\label{thm1c}
\end{equation}
Examining the latter integrand we see that
\begin{eqnarray}
 \nabla^{^{B_0}}_{_{(K'}}
 \nabla_{_{C')[B_0}} \bar{f}_{_{B_1 \ldots B_m]A'_2 \ldots A'_p
}}^{^{\;C'}}&=&
\frac{1}{m+1}( \nabla^{^{B_0}}_{_{(K'}}
 \nabla_{_{C')B_0}} \bar{f}_{_{B_1 \ldots B_mA'_2 \ldots A'_p}}^{^{\;C'}}
\nonumber \\
&&+ \sum_{i=1}^{i=m}(-1)^i \nabla^{^{B_0}}_{_{(K'}}
 \nabla_{_{C')B_i}} \bar{f}_{_{B_0B_1 \ldots \hat{B}_i \ldots B_m A'_2
\ldots A'_p}}^{^{\;C'}})
\label{thm1d}
\end{eqnarray}
Now by definition
\begin{displaymath}
\nabla^{^{B_0}}_{_{(K'}}\nabla_{_{C')B_i}}-\nabla_{_{B_i(K'}}
\nabla_{_{C')}}^{^{B_0}}=\Box_{_{K'C' \quad B_i}}^{^{ \qquad B_0}}.
\end{displaymath}
(For the definition and properties of the curvature operators see
\cite{pcf})
The second operator on the left above is
$ \frac{1}{2}(\nabla_{_{B_iK'}}\nabla_{_{C'}}^{^{B_0}}+\nabla_{_{B_iC'}}
\nabla_{_{K'}}^{^{B_0}})$ and the first part of this sum annihilates
$\bar{f}$
 by (c) of (\ref{hodge1}).
Equation (\ref{thm1c}) can now be rewritten
\begin{eqnarray*}
\frac{(p+2)}{2(p+1)} \parallel \nabla_{_{A'}}^{^{[B_0}}f^{^{B_1
\ldots B_m]A'
A'_2 \ldots A'_p}}\parallel^2&=& -\frac{1}{m+1} \int f^{^{B_1
\ldots B_m K'A'_2
\ldots A'_p}}(\Box_{_{K'C' \quad B_0}}
^{^{ \qquad B_0}} \bar{f}^{^{\; C'}}_{_{B_1 \ldots B_m A'_2
\ldots A'_p}})
 \\
-\frac{1}{(m+1)} \sum_{i=1}^{i=m}& (-1)^i & \int f^{^{B_1
\ldots B_m K'A'_2
\ldots A'_p}} (\Box_{_{K'C' \quad B_i}}^{^{ \qquad B_0}}
\bar{f}^{^{C'}}
_{_{B_0B_1 \ldots \hat{B}_i \ldots B_p A'_2 \ldots A'_p}}) \\
-\frac{1}{2(m+1)}\sum_{i=1}^{i=m} &(-1)^i & \int f^{^{B_1
\ldots B_m K'A'_2
\ldots A'_p}}(\nabla_{_{B_iC'}}\nabla^{^{B_0 }}_{_{K'}}
\bar{f}^{^{C'}}
_{_{B_0 B_1 \ldots \hat{B}_i \ldots B_m 'A'_2 \ldots A'_p}})
\end{eqnarray*}
Using the symmetry in the primed indices given by the latter part of
lemma~(\ref{lemma2}) we may commute the final $K'$ and $C'$ in the
final
integral above and, after another application of integration by parts,
this final
integral becomes
\begin{eqnarray*}
&- \int(\nabla _{_{B_iC'}}f^{^{B_1 \ldots B_m K' A'_2 \ldots A'_p}})
(\nabla
^{^{B_0C'}}  \bar{f}_{_{B_0 B_1 \ldots \hat{B}_i \ldots B_m K'A'_2
\ldots A'_p}})
& \\
&=(-1)^i\int (\nabla_{_{B_i}}^{^{C'}}f^{^{B_iB_1 \ldots \hat{B}_i
\ldots
B_m K'A'_2 \ldots A'_p}})( \nabla _{_{B_0C'}} \bar{f}^{^{B_0}}_{_{
\quad B_1 \ldots \hat{B} \ldots B_m  K' A'_2 \ldots A'_p }})&\\
&= (-1)^i \parallel \nabla_{_{B}}^{^{C'}}f^{^{BB_2 \ldots B_m A'_1
\ldots A'_p}}
\parallel ^2 &
\end{eqnarray*}
With this in mind we may now rearrange the above equation to obtain
\begin{equation}
\frac{(p+2)}{2(p+1)} \parallel \nabla_{_{A'}}^{^{[B_0}}f^{^{B_1
\ldots B_m]A'
A'_2 \ldots A'_p}}\parallel^2 + \frac{m}{2(m+1)}
\parallel \nabla_{_{B}}^{^{C'}}f^{^{BB_2 \ldots B_m A'_1 \ldots A'_p}}
\parallel ^2
= -\frac{1}{m+1} \int {\em R}(f) \label{thm1e}
\end{equation}
where ${\em R}(f)$ is given by
\begin{displaymath}
{\em R}(f)=f^{^{B_1 \ldots B_m K'A'_2
\ldots A'_p}}(\Box_{_{K'C' \quad B_0}}
^{^{ \qquad B_0}} \bar{f}^{^{\; C'}}_{_{B_1 \ldots B_m A'_2
\ldots A'_p}}+
\sum_{i=1}^{i=m}(-1)^i
\Box_{_{K'C' \quad B_i}}^{^{ \qquad B_0}} \bar{f}^{^{C'}}
_{_{B_0B_1 \ldots \hat{B}_i \ldots B_p A'_2 \ldots A'_p}})
\end{displaymath}
The evaluation of the action of the curvature operators is an elementary
exercise. We have
\begin{displaymath}
\Box_{_{K'C' \quad B_0}}^{^{ \qquad B_0}} \bar{f}^{^{\; C'}}_{_{B_1
\ldots B_m
 A'_2 \ldots A'_p}} =-2 \Lambda k(p+2) \bar{f}_{_{B_1 \ldots B_m K'
 A'_2 \ldots A'_p}}
\end{displaymath}
and the other operator gives
\begin{eqnarray*}
\Box_{_{K'C' \quad B_i}}^{^{ \qquad B_0}} \bar{f}
_{_{B_0 \ldots \hat{B}_i \ldots B_p}}^{^{ A'_1 \ldots A'_p}}&=&
-2 \Lambda e^{^{B_0}}_{_{\quad B_i}} \sum_{j=1}^{j=p} \eps_{_{(K'}}^{^{
\;\quad
A'_j}} \eps_{_{C')Q'}} \bar{f} ^{^{Q'A'_1 \ldots \hat{A}'_j \ldots A'_p}}
_{_{B_0 \ldots \hat{B}_i \ldots B_m}} \\
&=&2 \Lambda \sum_{j=1}^{j=p} \eps_{_{(K'}}^{^{\;\quad
A'_j}} \eps_{_{C')Q'}} \bar{f} ^{^{Q'A'_1 \ldots \hat{A}'_j \ldots A'_p}}
_{_{B_i B_1 \ldots \hat{B}_i \ldots B_m}}\\
&=&(-1)^{i-1}2 \Lambda \sum_{j=1}^{j=p} \eps_{_{(K'}}^{^{\;\quad
A'_j}} \eps_{_{C')Q'}} \bar{f} ^{^{Q'A'_1 \ldots \hat{A}'_j \ldots A'_p}}
_{_{ B_1 \ldots  B_m}}\\
 && \\
\Box_{_{K'C' \quad B_i}}^{^{\qquad B_0}} \bar{f}^{^{\; C'}}_{_{B_0 \ldots
\hat{B}_i \ldots B_m A'_2 \ldots A'_p}}&=& (-1)^i \Lambda (p+2)
 \bar{f} _{_{
B_1 \ldots B_m K'A'_2 \ldots A'_p}}
\end{eqnarray*}
so that ${\em R}(f)=-\Lambda (p+2)(2k-m)\bar{f} _{_{B_1 \ldots B_m K'A'_2
\ldots A'_p}}$. Substituting this into (\ref{thm1e}) completes the proof.
$\Box$


\begin{thebibliography}{99}

\bibitem{pcf} T.N.Bailey, M.G.Eastwood,  Complex Paraconformal
Manifolds -
 Their Differential Geometry and Twistor Theory, {\em Forum Math.},
{\bf 3}, (1991), 61-103.
\bibitem{hor} R.E.Horan, A Rigidity Theorem for Quaternionic-K\"{a}hler
 Manifolds, To appear in Differential Geometry and its Applications
\bibitem{sst} R.Penrose, W.Rindler, {\em Spinors and Space-Time}, Vol 1,
(C.U.P., Cambridge, 1984).
\bibitem{sal1} S.M.Salamon,  Quaternionic K\"{a}hler manifolds,
{\em Invent.
 Math.}, {\bf 67}, (1982), 143-171.
\bibitem{sal2} S.M.Salamon, Differential Geometry of Quaternionic
Manifolds,
 {\em Ann. Scient. Ec. Norm. Sup.}, {\bf 19}, (1986), 31-55.
\bibitem{wells} R.O.Wells, {\em Differential Analysis on Complex
Manifolds},
 (Springer-Verlag, Berlin, 1980).
\end{thebibliography}
\end{document}